\documentclass{amsart}
\usepackage{amsfonts}

\theoremstyle{plain}

\numberwithin{equation}{section}
\input{tcilatex}

\begin{document}
\title{Examples of Abel Grassmann's groupoids}
\author{}
\maketitle

\bigskip

\begin{center}
\bigskip

\textbf{Qaiser Mushtaq}

Department of Mathematics

Quaid-i-Azam University, Islamabad, Pakistan

\bigskip

\textbf{Madad Khan and Sameera Bano}

Department of Mathematics

COMSATS Institute of Information Technology, Abbottabad, Pakistan

\bigskip
\end{center}

Following are the examples of AG-groupoids(LA-semigroups), AG$^{\ast \ast }$%
-groupoids, Inverse AG-groupoids, Locally associative AG-groupoids with left
identity, Locally associative AG$^{\ast \ast }$-groupoids and AG-bands.

\begin{equation*}
\begin{tabular}{c|ccc}
$\cdot $ & $1$ & $2$ & $3$ \\ \hline
$1$ & $1$ & $1$ & $1$ \\ 
$2$ & $3$ & $3$ & $3$ \\ 
$3$ & $1$ & $1$ & $1$%
\end{tabular}%
\end{equation*}

\begin{equation*}
\begin{tabular}{c|ccc}
$\cdot $ & $1$ & $2$ & $3$ \\ \hline
$1$ & $2$ & $3$ & $2$ \\ 
$2$ & $2$ & $2$ & $2$ \\ 
$3$ & $2$ & $2$ & $2$%
\end{tabular}%
\end{equation*}

\begin{equation*}
\begin{tabular}{c|ccc}
$\cdot $ & $1$ & $2$ & $3$ \\ \hline
$1$ & $2$ & $2$ & $2$ \\ 
$2$ & $2$ & $2$ & $2$ \\ 
$3$ & $1$ & $1$ & $2$%
\end{tabular}%
\end{equation*}

\begin{equation*}
\begin{tabular}{c|ccc}
$\cdot $ & $1$ & $2$ & $3$ \\ \hline
$1$ & $3$ & $2$ & $3$ \\ 
$2$ & $2$ & $2$ & $2$ \\ 
$3$ & $2$ & $2$ & $2$%
\end{tabular}%
\end{equation*}

\begin{equation*}
\begin{tabular}{c|ccc}
$\cdot $ & $1$ & $2$ & $3$ \\ \hline
$1$ & $2$ & $3$ & $2$ \\ 
$2$ & $3$ & $3$ & $3$ \\ 
$3$ & $3$ & $3$ & $3$%
\end{tabular}%
\end{equation*}

\begin{equation*}
\begin{tabular}{c|ccc}
$\cdot $ & $1$ & $2$ & $3$ \\ \hline
$1$ & $1$ & $1$ & $1$ \\ 
$2$ & $3$ & $3$ & $1$ \\ 
$3$ & $1$ & $1$ & $1$%
\end{tabular}%
\end{equation*}

\begin{equation*}
\begin{tabular}{c|ccc}
$\cdot $ & $1$ & $2$ & $3$ \\ \hline
$1$ & $2$ & $2$ & $2$ \\ 
$2$ & $3$ & $3$ & $3$ \\ 
$3$ & $3$ & $3$ & $3$%
\end{tabular}%
\end{equation*}

\begin{equation*}
\begin{tabular}{c|ccc}
$\cdot $ & $1$ & $2$ & $3$ \\ \hline
$1$ & $3$ & $3$ & $3$ \\ 
$2$ & $1$ & $3$ & $1$ \\ 
$3$ & $3$ & $3$ & $3$%
\end{tabular}%
\end{equation*}

\begin{equation*}
\begin{tabular}{c|ccc}
$\cdot $ & $1$ & $2$ & $3$ \\ \hline
$1$ & $1$ & $1$ & $1$ \\ 
$2$ & $1$ & $3$ & $3$ \\ 
$3$ & $1$ & $1$ & $1$%
\end{tabular}%
\end{equation*}

\begin{equation*}
\begin{tabular}{c|ccc}
$\cdot $ & $1$ & $2$ & $3$ \\ \hline
$1$ & $3$ & $3$ & $2$ \\ 
$2$ & $3$ & $3$ & $3$ \\ 
$3$ & $3$ & $3$ & $3$%
\end{tabular}%
\end{equation*}

\begin{equation*}
\begin{tabular}{c|ccc}
$\cdot $ & $1$ & $2$ & $3$ \\ \hline
$1$ & $3$ & $2$ & $2$ \\ 
$2$ & $3$ & $3$ & $3$ \\ 
$3$ & $3$ & $3$ & $3$%
\end{tabular}%
\end{equation*}

\begin{equation*}
\begin{tabular}{c|ccc}
$\cdot $ & $1$ & $2$ & $3$ \\ \hline
$1$ & $1$ & $1$ & $1$ \\ 
$2$ & $1$ & $1$ & $1$ \\ 
$3$ & $1$ & $2$ & $1$%
\end{tabular}%
\end{equation*}

\begin{equation*}
\begin{tabular}{c|ccc}
$\cdot $ & $1$ & $2$ & $3$ \\ \hline
$1$ & $3$ & $2$ & $3$ \\ 
$2$ & $3$ & $3$ & $3$ \\ 
$3$ & $3$ & $3$ & $3$%
\end{tabular}%
\end{equation*}

\begin{equation*}
\begin{tabular}{c|ccc}
$\cdot $ & $1$ & $2$ & $3$ \\ \hline
$1$ & $2$ & $2$ & $3$ \\ 
$2$ & $3$ & $3$ & $3$ \\ 
$3$ & $3$ & $3$ & $3$%
\end{tabular}%
\end{equation*}

\begin{equation*}
\begin{tabular}{c|ccc}
$\cdot $ & $1$ & $2$ & $3$ \\ \hline
$1$ & $2$ & $2$ & $3$ \\ 
$2$ & $2$ & $2$ & $2$ \\ 
$3$ & $2$ & $2$ & $2$%
\end{tabular}%
\end{equation*}

\begin{equation*}
\begin{tabular}{c|ccc}
$\cdot $ & $1$ & $2$ & $3$ \\ \hline
$1$ & $3$ & $3$ & $3$ \\ 
$2$ & $3$ & $1$ & $1$ \\ 
$3$ & $3$ & $3$ & $3$%
\end{tabular}%
\end{equation*}

\begin{equation*}
\begin{tabular}{c|ccc}
$\cdot $ & $1$ & $2$ & $3$ \\ \hline
$1$ & $3$ & $3$ & $2$ \\ 
$2$ & $2$ & $2$ & $2$ \\ 
$3$ & $2$ & $2$ & $2$%
\end{tabular}%
\end{equation*}

\begin{equation*}
\begin{tabular}{c|ccc}
$\cdot $ & $1$ & $2$ & $3$ \\ \hline
$1$ & $2$ & $2$ & $2$ \\ 
$2$ & $2$ & $2$ & $2$ \\ 
$3$ & $1$ & $2$ & $1$%
\end{tabular}%
\end{equation*}

\begin{equation*}
\begin{tabular}{c|ccc}
$\cdot $ & $1$ & $2$ & $3$ \\ \hline
$1$ & $1$ & $1$ & $1$ \\ 
$2$ & $1$ & $1$ & $1$ \\ 
$3$ & $2$ & $1$ & $2$%
\end{tabular}%
\end{equation*}

\begin{equation*}
\begin{tabular}{c|ccc}
$\cdot $ & $1$ & $2$ & $3$ \\ \hline
$1$ & $3$ & $3$ & $3$ \\ 
$2$ & $1$ & $1$ & $3$ \\ 
$3$ & $3$ & $3$ & $3$%
\end{tabular}%
\end{equation*}

\begin{equation*}
\begin{tabular}{c|ccc}
$\cdot $ & $1$ & $2$ & $3$ \\ \hline
$1$ & $3$ & $3$ & $3$ \\ 
$2$ & $2$ & $2$ & $2$ \\ 
$3$ & $2$ & $2$ & $2$%
\end{tabular}%
\end{equation*}

\begin{equation*}
\begin{tabular}{c|ccc}
$\cdot $ & $1$ & $2$ & $3$ \\ \hline
$1$ & $2$ & $2$ & $2$ \\ 
$2$ & $2$ & $2$ & $2$ \\ 
$3$ & $1$ & $1$ & $1$%
\end{tabular}%
\end{equation*}

\begin{equation*}
\begin{tabular}{c|ccc}
$\cdot $ & $1$ & $2$ & $3$ \\ \hline
$1$ & $2$ & $2$ & $2$ \\ 
$2$ & $2$ & $2$ & $2$ \\ 
$3$ & $2$ & $1$ & $1$%
\end{tabular}%
\end{equation*}

\begin{equation*}
\begin{tabular}{c|ccc}
$\cdot $ & $1$ & $2$ & $3$ \\ \hline
$1$ & $1$ & $1$ & $1$ \\ 
$2$ & $1$ & $1$ & $1$ \\ 
$3$ & $2$ & $2$ & $2$%
\end{tabular}%
\end{equation*}

\begin{equation*}
\begin{tabular}{c|ccc}
$\cdot $ & $1$ & $2$ & $3$ \\ \hline
$1$ & $3$ & $3$ & $3$ \\ 
$2$ & $1$ & $1$ & $1$ \\ 
$3$ & $3$ & $3$ & $3$%
\end{tabular}%
\end{equation*}

\begin{equation*}
\begin{tabular}{c|ccc}
$\cdot $ & $1$ & $2$ & $3$ \\ \hline
$1$ & $1$ & $1$ & $1$ \\ 
$2$ & $3$ & $1$ & $1$ \\ 
$3$ & $1$ & $1$ & $1$%
\end{tabular}%
\end{equation*}

\begin{equation*}
\begin{tabular}{c|ccc}
$\cdot $ & $1$ & $2$ & $3$ \\ \hline
$1$ & $3$ & $3$ & $3$ \\ 
$2$ & $3$ & $3$ & $3$ \\ 
$3$ & $1$ & $3$ & $3$%
\end{tabular}%
\end{equation*}

\begin{equation*}
\begin{tabular}{c|ccc}
$\cdot $ & $1$ & $2$ & $3$ \\ \hline
$1$ & $1$ & $1$ & $1$ \\ 
$2$ & $1$ & $1$ & $1$ \\ 
$3$ & $1$ & $2$ & $2$%
\end{tabular}%
\end{equation*}

\begin{equation*}
\begin{tabular}{c|ccc}
$\cdot $ & $1$ & $2$ & $3$ \\ \hline
$1$ & $1$ & $1$ & $1$ \\ 
$2$ & $1$ & $1$ & $1$ \\ 
$3$ & $2$ & $2$ & $1$%
\end{tabular}%
\end{equation*}

\begin{equation*}
\begin{tabular}{c|ccc}
$\cdot $ & $1$ & $2$ & $3$ \\ \hline
$1$ & $2$ & $3$ & $3$ \\ 
$2$ & $2$ & $2$ & $2$ \\ 
$3$ & $2$ & $2$ & $2$%
\end{tabular}%
\end{equation*}

\begin{equation*}
\begin{tabular}{c|ccc}
$\cdot $ & $1$ & $2$ & $3$ \\ \hline
$1$ & $1$ & $1$ & $1$ \\ 
$2$ & $3$ & $1$ & $3$ \\ 
$3$ & $1$ & $1$ & $1$%
\end{tabular}%
\end{equation*}

\begin{equation*}
\begin{tabular}{c|ccc}
$\cdot $ & $1$ & $2$ & $3$ \\ \hline
$1$ & $3$ & $3$ & $3$ \\ 
$2$ & $1$ & $3$ & $3$ \\ 
$3$ & $3$ & $3$ & $3$%
\end{tabular}%
\end{equation*}

\begin{equation*}
\begin{tabular}{c|ccc}
$\cdot $ & $1$ & $2$ & $3$ \\ \hline
$1$ & $1$ & $1$ & $1$ \\ 
$2$ & $1$ & $1$ & $1$ \\ 
$3$ & $2$ & $1$ & $1$%
\end{tabular}%
\end{equation*}

\begin{equation*}
\begin{tabular}{c|ccc}
$\cdot $ & $1$ & $2$ & $3$ \\ \hline
$1$ & $1$ & $1$ & $1$ \\ 
$2$ & $3$ & $1$ & $2$ \\ 
$3$ & $1$ & $1$ & $1$%
\end{tabular}%
\end{equation*}

\begin{equation*}
\begin{tabular}{c|ccc}
$\cdot $ & $1$ & $2$ & $3$ \\ \hline
$1$ & $2$ & $3$ & $1$ \\ 
$2$ & $1$ & $2$ & $3$ \\ 
$3$ & $3$ & $1$ & $2$%
\end{tabular}%
\end{equation*}

\begin{equation*}
\begin{tabular}{c|ccc}
$\cdot $ & $1$ & $2$ & $3$ \\ \hline
$1$ & $2$ & $2$ & $3$ \\ 
$2$ & $2$ & $2$ & $2$ \\ 
$3$ & $1$ & $2$ & $2$%
\end{tabular}%
\end{equation*}

\begin{equation*}
\begin{tabular}{c|ccc}
$\cdot $ & $1$ & $2$ & $3$ \\ \hline
$1$ & $3$ & $1$ & $2$ \\ 
$2$ & $2$ & $3$ & $1$ \\ 
$3$ & $1$ & $2$ & $3$%
\end{tabular}%
\end{equation*}

\begin{equation*}
\begin{tabular}{c|ccc}
$\cdot $ & $1$ & $2$ & $3$ \\ \hline
$1$ & $1$ & $1$ & $1$ \\ 
$2$ & $1$ & $1$ & $3$ \\ 
$3$ & $1$ & $2$ & $1$%
\end{tabular}%
\end{equation*}

\begin{equation*}
\begin{tabular}{c|cccc}
$\cdot $ & $1$ & $2$ & $3$ & $4$ \\ \hline
$1$ & $4$ & $4$ & $2$ & $2$ \\ 
$2$ & $4$ & $4$ & $4$ & $4$ \\ 
$3$ & $4$ & $4$ & $2$ & $4$ \\ 
$4$ & $4$ & $4$ & $4$ & $4$%
\end{tabular}%
\end{equation*}

\begin{equation*}
\begin{tabular}{c|cccc}
$\cdot $ & $1$ & $2$ & $3$ & $4$ \\ \hline
$1$ & $1$ & $1$ & $1$ & $1$ \\ 
$2$ & $1$ & $1$ & $1$ & $1$ \\ 
$3$ & $2$ & $2$ & $2$ & $2$ \\ 
$4$ & $2$ & $2$ & $1$ & $1$%
\end{tabular}%
\end{equation*}

\begin{equation*}
\begin{tabular}{c|cccc}
$\cdot $ & $1$ & $2$ & $3$ & $4$ \\ \hline
$1$ & $1$ & $1$ & $1$ & $1$ \\ 
$2$ & $4$ & $1$ & $4$ & $1$ \\ 
$3$ & $1$ & $4$ & $4$ & $1$ \\ 
$4$ & $1$ & $1$ & $1$ & $1$%
\end{tabular}%
\end{equation*}

\begin{equation*}
\begin{tabular}{c|cccc}
$\cdot $ & $1$ & $2$ & $3$ & $4$ \\ \hline
$1$ & $1$ & $1$ & $1$ & $1$ \\ 
$2$ & $1$ & $1$ & $1$ & $1$ \\ 
$3$ & $2$ & $1$ & $2$ & $2$ \\ 
$4$ & $1$ & $1$ & $1$ & $1$%
\end{tabular}%
\end{equation*}

\begin{equation*}
\begin{tabular}{c|cccc}
$\cdot $ & $1$ & $2$ & $3$ & $4$ \\ \hline
$1$ & $4$ & $3$ & $3$ & $3$ \\ 
$2$ & $3$ & $3$ & $4$ & $4$ \\ 
$3$ & $3$ & $3$ & $3$ & $3$ \\ 
$4$ & $3$ & $3$ & $3$ & $3$%
\end{tabular}%
\end{equation*}

\begin{equation*}
\begin{tabular}{c|cccc}
$\cdot $ & $1$ & $2$ & $3$ & $4$ \\ \hline
$1$ & $1$ & $1$ & $1$ & $1$ \\ 
$2$ & $1$ & $2$ & $3$ & $4$ \\ 
$3$ & $1$ & $4$ & $2$ & $3$ \\ 
$4$ & $1$ & $3$ & $4$ & $2$%
\end{tabular}%
\end{equation*}

\begin{equation*}
\begin{tabular}{c|cccc}
$\cdot $ & $1$ & $2$ & $3$ & $4$ \\ \hline
$1$ & $1$ & $1$ & $1$ & $1$ \\ 
$2$ & $1$ & $2$ & $3$ & $4$ \\ 
$3$ & $1$ & $3$ & $4$ & $2$ \\ 
$4$ & $1$ & $4$ & $2$ & $3$%
\end{tabular}%
\end{equation*}

\begin{equation*}
\begin{tabular}{c|cccc}
$\cdot $ & $1$ & $2$ & $3$ & $4$ \\ \hline
$1$ & $1$ & $1$ & $1$ & $1$ \\ 
$2$ & $1$ & $4$ & $2$ & $3$ \\ 
$3$ & $1$ & $2$ & $3$ & $4$ \\ 
$4$ & $1$ & $3$ & $4$ & $2$%
\end{tabular}%
\end{equation*}

\begin{equation*}
\begin{tabular}{c|cccc}
$\cdot $ & $1$ & $2$ & $3$ & $4$ \\ \hline
$1$ & $1$ & $1$ & $1$ & $1$ \\ 
$2$ & $1$ & $3$ & $4$ & $2$ \\ 
$3$ & $1$ & $4$ & $2$ & $3$ \\ 
$4$ & $1$ & $2$ & $3$ & $4$%
\end{tabular}%
\end{equation*}

\bigskip

\begin{equation*}
\begin{tabular}{c|cccc}
$\cdot $ & $a$ & $b$ & $c$ & $d$ \\ \hline
$a$ & $d$ & $a$ & $b$ & $c$ \\ 
$b$ & $c$ & $d$ & $a$ & $b$ \\ 
$c$ & $b$ & $c$ & $d$ & $a$ \\ 
$d$ & $a$ & $b$ & $c$ & $d$%
\end{tabular}%
\end{equation*}

\begin{equation*}
\begin{tabular}{c|cccc}
$\cdot $ & $a$ & $b$ & $c$ & $d$ \\ \hline
$a$ & $b$ & $c$ & $d$ & $a$ \\ 
$b$ & $a$ & $b$ & $c$ & $d$ \\ 
$c$ & $d$ & $a$ & $b$ & $c$ \\ 
$d$ & $c$ & $d$ & $a$ & $b$%
\end{tabular}%
\end{equation*}

\begin{equation*}
\begin{tabular}{c|cccc}
$\cdot $ & $1$ & $2$ & $3$ & $4$ \\ \hline
$1$ & $1$ & $1$ & $1$ & $1$ \\ 
$2$ & $1$ & $1$ & $1$ & $1$ \\ 
$3$ & $2$ & $1$ & $2$ & $1$ \\ 
$4$ & $1$ & $2$ & $2$ & $2$%
\end{tabular}%
\end{equation*}%
\begin{equation*}
\begin{tabular}{c|cccc}
$\cdot $ & $1$ & $2$ & $3$ & $4$ \\ \hline
$1$ & $2$ & $2$ & $2$ & $2$ \\ 
$2$ & $2$ & $2$ & $2$ & $2$ \\ 
$3$ & $1$ & $4$ & $4$ & $2$ \\ 
$4$ & $2$ & $2$ & $2$ & $2$%
\end{tabular}%
\end{equation*}

\begin{equation*}
\begin{tabular}{c|cccc}
$\cdot $ & $1$ & $2$ & $3$ & $4$ \\ \hline
$1$ & $4$ & $3$ & $4$ & $4$ \\ 
$2$ & $1$ & $4$ & $4$ & $1$ \\ 
$3$ & $4$ & $4$ & $4$ & $4$ \\ 
$4$ & $4$ & $4$ & $4$ & $4$%
\end{tabular}%
\end{equation*}

\begin{equation*}
\begin{tabular}{c|cccc}
$\cdot $ & $1$ & $2$ & $3$ & $4$ \\ \hline
$1$ & $3$ & $2$ & $3$ & $2$ \\ 
$2$ & $2$ & $2$ & $2$ & $2$ \\ 
$3$ & $2$ & $2$ & $2$ & $2$ \\ 
$4$ & $2$ & $2$ & $3$ & $2$%
\end{tabular}%
\end{equation*}

\begin{equation*}
\begin{tabular}{c|cccc}
$\cdot $ & $1$ & $2$ & $3$ & $4$ \\ \hline
$1$ & $1$ & $1$ & $1$ & $1$ \\ 
$2$ & $1$ & $1$ & $4$ & $1$ \\ 
$3$ & $1$ & $4$ & $1$ & $2$ \\ 
$4$ & $1$ & $1$ & $1$ & $1$%
\end{tabular}%
\end{equation*}

\begin{equation*}
\begin{tabular}{c|cccc}
$\cdot $ & $1$ & $2$ & $3$ & $4$ \\ \hline
$1$ & $1$ & $1$ & $1$ & $1$ \\ 
$2$ & $3$ & $3$ & $4$ & $4$ \\ 
$3$ & $1$ & $4$ & $1$ & $1$ \\ 
$4$ & $1$ & $1$ & $1$ & $1$%
\end{tabular}%
\end{equation*}

\begin{equation*}
\begin{tabular}{c|cccc}
$\cdot $ & $1$ & $2$ & $3$ & $4$ \\ \hline
$1$ & $2$ & $4$ & $3$ & $2$ \\ 
$2$ & $4$ & $4$ & $4$ & $4$ \\ 
$3$ & $4$ & $4$ & $4$ & $4$ \\ 
$4$ & $4$ & $4$ & $4$ & $4$%
\end{tabular}%
\end{equation*}

\begin{equation*}
\begin{tabular}{c|cccc}
$\cdot $ & $1$ & $2$ & $3$ & $4$ \\ \hline
$1$ & $2$ & $4$ & $4$ & $4$ \\ 
$2$ & $2$ & $2$ & $2$ & $2$ \\ 
$3$ & $4$ & $2$ & $4$ & $4$ \\ 
$4$ & $2$ & $2$ & $2$ & $2$%
\end{tabular}%
\end{equation*}

\begin{equation*}
\begin{tabular}{c|cccc}
$\cdot $ & $1$ & $2$ & $3$ & $4$ \\ \hline
$1$ & $4$ & $4$ & $4$ & $4$ \\ 
$2$ & $1$ & $4$ & $4$ & $4$ \\ 
$3$ & $1$ & $4$ & $4$ & $4$ \\ 
$4$ & $4$ & $4$ & $4$ & $4$%
\end{tabular}%
\end{equation*}

\begin{equation*}
\begin{tabular}{c|cccc}
$\cdot $ & $1$ & $2$ & $3$ & $4$ \\ \hline
$1$ & $4$ & $4$ & $4$ & $4$ \\ 
$2$ & $1$ & $1$ & $1$ & $1$ \\ 
$3$ & $1$ & $4$ & $1$ & $4$ \\ 
$4$ & $4$ & $4$ & $4$ & $4$%
\end{tabular}%
\end{equation*}

\begin{equation*}
\begin{tabular}{c|cccc}
$\cdot $ & $1$ & $2$ & $3$ & $4$ \\ \hline
$1$ & $4$ & $4$ & $4$ & $4$ \\ 
$2$ & $3$ & $1$ & $3$ & $1$ \\ 
$3$ & $4$ & $1$ & $4$ & $4$ \\ 
$4$ & $4$ & $4$ & $4$ & $4$%
\end{tabular}%
\end{equation*}

\begin{equation*}
\begin{tabular}{c|cccc}
$\cdot $ & $a$ & $b$ & $c$ & $d$ \\ \hline
$a$ & $b$ & $b$ & $c$ & $c$ \\ 
$b$ & $b$ & $b$ & $b$ & $b$ \\ 
$c$ & $b$ & $b$ & $b$ & $b$ \\ 
$d$ & $b$ & $b$ & $b$ & $b$%
\end{tabular}%
\end{equation*}

\begin{equation*}
\begin{tabular}{c|cccc}
$\cdot $ & $a$ & $b$ & $c$ & $d$ \\ \hline
$a$ & $b$ & $b$ & $b$ & $c$ \\ 
$b$ & $b$ & $b$ & $b$ & $b$ \\ 
$c$ & $b$ & $b$ & $b$ & $b$ \\ 
$d$ & $b$ & $c$ & $b$ & $b$%
\end{tabular}%
\end{equation*}

\begin{equation*}
\begin{tabular}{c|cccc}
$\cdot $ & $a$ & $b$ & $c$ & $d$ \\ \hline
$a$ & $b$ & $b$ & $d$ & $c$ \\ 
$b$ & $b$ & $b$ & $b$ & $b$ \\ 
$c$ & $b$ & $b$ & $b$ & $b$ \\ 
$d$ & $b$ & $b$ & $b$ & $b$%
\end{tabular}%
\end{equation*}

\begin{equation*}
\begin{tabular}{c|cccc}
$\cdot $ & $a$ & $b$ & $c$ & $d$ \\ \hline
$a$ & $a$ & $a$ & $a$ & $a$ \\ 
$b$ & $a$ & $a$ & $a$ & $a$ \\ 
$c$ & $a$ & $a$ & $a$ & $a$ \\ 
$d$ & $a$ & $a$ & $c$ & $a$%
\end{tabular}%
\end{equation*}

\begin{equation*}
\begin{tabular}{c|cccc}
$\cdot $ & $a$ & $b$ & $c$ & $d$ \\ \hline
$a$ & $a$ & $a$ & $a$ & $a$ \\ 
$b$ & $a$ & $a$ & $c$ & $a$ \\ 
$c$ & $a$ & $a$ & $a$ & $a$ \\ 
$d$ & $a$ & $a$ & $a$ & $a$%
\end{tabular}%
\end{equation*}

\begin{equation*}
\begin{tabular}{c|cccc}
$\cdot $ & $a$ & $b$ & $c$ & $d$ \\ \hline
$a$ & $a$ & $a$ & $a$ & $a$ \\ 
$b$ & $c$ & $a$ & $a$ & $a$ \\ 
$c$ & $a$ & $a$ & $a$ & $a$ \\ 
$d$ & $a$ & $a$ & $a$ & $a$%
\end{tabular}%
\end{equation*}

\begin{equation*}
\begin{tabular}{c|cccc}
$\cdot $ & $a$ & $b$ & $c$ & $d$ \\ \hline
$a$ & $a$ & $a$ & $a$ & $a$ \\ 
$b$ & $a$ & $a$ & $a$ & $a$ \\ 
$c$ & $a$ & $a$ & $a$ & $a$ \\ 
$d$ & $c$ & $a$ & $a$ & $a$%
\end{tabular}%
\end{equation*}

\begin{equation*}
\begin{tabular}{c|cccc}
$\cdot $ & $a$ & $b$ & $c$ & $d$ \\ \hline
$a$ & $b$ & $b$ & $b$ & $b$ \\ 
$b$ & $b$ & $b$ & $b$ & $b$ \\ 
$c$ & $b$ & $b$ & $b$ & $b$ \\ 
$d$ & $b$ & $c$ & $b$ & $b$%
\end{tabular}%
\end{equation*}

\begin{equation*}
\begin{tabular}{c|cccc}
$\cdot $ & $a$ & $b$ & $c$ & $d$ \\ \hline
$a$ & $b$ & $c$ & $b$ & $b$ \\ 
$b$ & $b$ & $b$ & $b$ & $b$ \\ 
$c$ & $b$ & $b$ & $b$ & $b$ \\ 
$d$ & $b$ & $b$ & $b$ & $b$%
\end{tabular}%
\end{equation*}

\begin{equation*}
\begin{tabular}{c|cccc}
$\cdot $ & $a$ & $b$ & $c$ & $d$ \\ \hline
$a$ & $b$ & $b$ & $b$ & $b$ \\ 
$b$ & $b$ & $b$ & $b$ & $b$ \\ 
$c$ & $b$ & $b$ & $b$ & $b$ \\ 
$d$ & $b$ & $b$ & $c$ & $b$%
\end{tabular}%
\end{equation*}

\begin{equation*}
\begin{tabular}{c|cccc}
$\cdot $ & $a$ & $b$ & $c$ & $d$ \\ \hline
$a$ & $b$ & $b$ & $c$ & $b$ \\ 
$b$ & $b$ & $b$ & $b$ & $b$ \\ 
$c$ & $b$ & $b$ & $b$ & $b$ \\ 
$d$ & $b$ & $b$ & $b$ & $b$%
\end{tabular}%
\end{equation*}

\begin{equation*}
\begin{tabular}{c|ccccc}
$\cdot $ & $1$ & $2$ & $3$ & $4$ & $5$ \\ \hline
$1$ & $1$ & $1$ & $1$ & $1$ & $1$ \\ 
$2$ & $1$ & $2$ & $2$ & $2$ & $2$ \\ 
$3$ & $1$ & $2$ & $4$ & $5$ & $3$ \\ 
$4$ & $1$ & $2$ & $3$ & $4$ & $5$ \\ 
$5$ & $1$ & $2$ & $5$ & $3$ & $4$%
\end{tabular}%
\end{equation*}

\begin{equation*}
\begin{tabular}{c|ccccc}
$\cdot $ & $1$ & $2$ & $3$ & $4$ & $5$ \\ \hline
$1$ & $1$ & $1$ & $1$ & $1$ & $1$ \\ 
$2$ & $1$ & $2$ & $2$ & $2$ & $2$ \\ 
$3$ & $1$ & $2$ & $4$ & $5$ & $3$ \\ 
$4$ & $1$ & $2$ & $5$ & $3$ & $4$ \\ 
$5$ & $1$ & $2$ & $3$ & $4$ & $5$%
\end{tabular}%
\end{equation*}

\begin{equation*}
\begin{tabular}{c|ccccc}
$\cdot $ & $1$ & $2$ & $3$ & $4$ & $5$ \\ \hline
$1$ & $2$ & $1$ & $1$ & $1$ & $1$ \\ 
$2$ & $1$ & $2$ & $2$ & $2$ & $2$ \\ 
$3$ & $1$ & $2$ & $3$ & $4$ & $5$ \\ 
$4$ & $1$ & $2$ & $5$ & $3$ & $4$ \\ 
$5$ & $1$ & $2$ & $4$ & $5$ & $3$%
\end{tabular}%
\end{equation*}

\begin{equation*}
\begin{tabular}{c|ccccc}
$\cdot $ & $1$ & $2$ & $3$ & $4$ & $5$ \\ \hline
$1$ & $1$ & $2$ & $3$ & $4$ & $5$ \\ 
$2$ & $5$ & $1$ & $2$ & $3$ & $4$ \\ 
$3$ & $4$ & $5$ & $1$ & $2$ & $3$ \\ 
$4$ & $3$ & $4$ & $5$ & $1$ & $2$ \\ 
$5$ & $2$ & $3$ & $4$ & $5$ & $1$%
\end{tabular}%
\end{equation*}

\begin{equation*}
\begin{tabular}{c|ccccc}
$\cdot $ & $1$ & $2$ & $3$ & $4$ & $5$ \\ \hline
$1$ & $5$ & $1$ & $2$ & $3$ & $4$ \\ 
$2$ & $4$ & $5$ & $1$ & $2$ & $3$ \\ 
$3$ & $3$ & $4$ & $5$ & $1$ & $2$ \\ 
$4$ & $2$ & $3$ & $4$ & $5$ & $1$ \\ 
$5$ & $1$ & $2$ & $3$ & $4$ & $5$%
\end{tabular}%
\end{equation*}

\begin{equation*}
\begin{tabular}{c|ccccc}
$\cdot $ & $1$ & $2$ & $3$ & $4$ & $5$ \\ \hline
$1$ & $4$ & $5$ & $1$ & $2$ & $3$ \\ 
$2$ & $3$ & $4$ & $5$ & $1$ & $2$ \\ 
$3$ & $2$ & $3$ & $4$ & $5$ & $1$ \\ 
$4$ & $1$ & $2$ & $3$ & $4$ & $5$ \\ 
$5$ & $5$ & $1$ & $2$ & $3$ & $4$%
\end{tabular}%
\end{equation*}

\begin{equation*}
\begin{tabular}{c|ccccc}
$\cdot $ & $1$ & $2$ & $3$ & $4$ & $5$ \\ \hline
$1$ & $3$ & $4$ & $5$ & $1$ & $2$ \\ 
$2$ & $2$ & $3$ & $4$ & $5$ & $1$ \\ 
$3$ & $1$ & $2$ & $3$ & $4$ & $5$ \\ 
$4$ & $5$ & $1$ & $2$ & $3$ & $4$ \\ 
$5$ & $4$ & $5$ & $1$ & $2$ & $3$%
\end{tabular}%
\end{equation*}

\begin{equation*}
\begin{tabular}{c|ccccc}
$\cdot $ & $1$ & $2$ & $3$ & $4$ & $5$ \\ \hline
$1$ & $2$ & $3$ & $4$ & $5$ & $1$ \\ 
$2$ & $1$ & $2$ & $3$ & $4$ & $5$ \\ 
$3$ & $5$ & $1$ & $2$ & $3$ & $4$ \\ 
$4$ & $4$ & $5$ & $1$ & $2$ & $3$ \\ 
$5$ & $3$ & $4$ & $5$ & $1$ & $2$%
\end{tabular}%
\end{equation*}

\begin{equation*}
\begin{tabular}{c|ccccc}
$\cdot $ & $1$ & $2$ & $3$ & $4$ & $5$ \\ \hline
$1$ & $1$ & $1$ & $1$ & $1$ & $1$ \\ 
$2$ & $1$ & $1$ & $1$ & $1$ & $1$ \\ 
$3$ & $1$ & $1$ & $4$ & $5$ & $3$ \\ 
$4$ & $1$ & $1$ & $3$ & $4$ & $5$ \\ 
$5$ & $1$ & $1$ & $5$ & $3$ & $4$%
\end{tabular}%
\end{equation*}

\begin{equation*}
\begin{tabular}{c|ccccc}
$\cdot $ & $1$ & $2$ & $3$ & $4$ & $5$ \\ \hline
$1$ & $2$ & $1$ & $1$ & $1$ & $1$ \\ 
$2$ & $1$ & $2$ & $2$ & $2$ & $2$ \\ 
$3$ & $1$ & $2$ & $4$ & $5$ & $3$ \\ 
$4$ & $1$ & $2$ & $3$ & $4$ & $5$ \\ 
$5$ & $1$ & $2$ & $5$ & $3$ & $4$%
\end{tabular}%
\end{equation*}

\begin{equation*}
\begin{tabular}{c|cccccc}
$\cdot $ & $1$ & $2$ & $3$ & $4$ & $5$ & $6$ \\ \hline
$1$ & $4$ & $4$ & $4$ & $4$ & $4$ & $6$ \\ 
$2$ & $4$ & $4$ & $4$ & $4$ & $4$ & $4$ \\ 
$3$ & $4$ & $4$ & $4$ & $4$ & $4$ & $4$ \\ 
$4$ & $4$ & $4$ & $4$ & $4$ & $4$ & $4$ \\ 
$5$ & $4$ & $4$ & $4$ & $4$ & $4$ & $4$ \\ 
$6$ & $4$ & $4$ & $4$ & $4$ & $4$ & $4$%
\end{tabular}%
\end{equation*}

\begin{equation*}
\begin{tabular}{c|cccccc}
$\cdot $ & $1$ & $2$ & $3$ & $4$ & $5$ & $6$ \\ \hline
$1$ & $4$ & $4$ & $4$ & $4$ & $4$ & $4$ \\ 
$2$ & $4$ & $4$ & $4$ & $4$ & $4$ & $4$ \\ 
$3$ & $4$ & $4$ & $4$ & $4$ & $4$ & $4$ \\ 
$4$ & $4$ & $4$ & $4$ & $4$ & $4$ & $4$ \\ 
$5$ & $4$ & $2$ & $4$ & $4$ & $4$ & $4$ \\ 
$6$ & $4$ & $4$ & $4$ & $4$ & $4$ & $4$%
\end{tabular}%
\end{equation*}

\begin{equation*}
\begin{tabular}{c|cccccc}
$\cdot $ & $1$ & $2$ & $3$ & $4$ & $5$ & $6$ \\ \hline
$1$ & $4$ & $4$ & $4$ & $4$ & $4$ & $4$ \\ 
$2$ & $4$ & $4$ & $4$ & $4$ & $5$ & $4$ \\ 
$3$ & $4$ & $4$ & $4$ & $4$ & $4$ & $4$ \\ 
$4$ & $4$ & $4$ & $4$ & $4$ & $4$ & $4$ \\ 
$5$ & $4$ & $4$ & $4$ & $4$ & $4$ & $4$ \\ 
$6$ & $4$ & $4$ & $4$ & $4$ & $4$ & $4$%
\end{tabular}%
\end{equation*}

\begin{equation*}
\begin{tabular}{c|ccccccc}
$\cdot $ & $1$ & $2$ & $3$ & $4$ & $5$ & $6$ & $7$ \\ \hline
$1$ & $4$ & $4$ & $4$ & $4$ & $4$ & $4$ & $4$ \\ 
$2$ & $4$ & $4$ & $4$ & $4$ & $4$ & $4$ & $4$ \\ 
$3$ & $4$ & $4$ & $4$ & $4$ & $4$ & $4$ & $4$ \\ 
$4$ & $4$ & $4$ & $4$ & $4$ & $4$ & $4$ & $4$ \\ 
$5$ & $4$ & $4$ & $4$ & $4$ & $4$ & $4$ & $4$ \\ 
$6$ & $4$ & $4$ & $4$ & $4$ & $4$ & $4$ & $4$ \\ 
$7$ & $1$ & $2$ & $3$ & $4$ & $5$ & $6$ & $4$%
\end{tabular}%
\end{equation*}

\end{document}